\pdfoutput=1
\documentclass[twoside]{article}
\makeatletter
\def\ps@myheadings{%
     \let\@oddfoot\@empty\let\@evenfoot\@empty
     \def\@evenhead{\thepage\hfil\MakeUppercase{\footnotesize\leftmark}\hfil}%
     \def\@oddhead{{\hfil\MakeUppercase{\footnotesize\rightmark}}\hfil\thepage}%
     \let\@mkboth\@gobbletwo
     \let\sectionmark\@gobble
     \let\subsectionmark\@gobble
     }
\makeatother
\usepackage{amsmath}
\usepackage{amsthm}
\usepackage{amscd}
\usepackage{amssymb}
\usepackage{amsfonts}
\usepackage{amscd}
\usepackage{amsbsy}
\usepackage{mathrsfs}
\usepackage{comment}

\usepackage{epsf}
\usepackage{epsfig}
\usepackage{psfig}
\usepackage{graphicx}

\newtheorem{theorem}{Theorem}[section]
\newtheorem{lemma}[theorem]{Lemma}
\newtheorem{prop}[theorem]{Proposition}

\theoremstyle{definition}

\newtheorem{quest}[theorem]{Question}

\theoremstyle{remark}
\newtheorem{remark}[theorem]{Remark}

\numberwithin{equation}{section}

\newcommand{\al}{\alpha}
\newcommand{\be}{\beta}
\newcommand{\Ga}{\Gamma}
\newcommand{\ga}{\gamma}
\newcommand{\sa}{\sigma}
\newcommand{\Sa}{\Sigma}

\newcommand{\La}{\Lambda}

\newcommand{\Om}{\Omega}

\newcommand{\scr}{{\mathcal R}}

\newcommand{\p}{\partial}

\newcommand{\ra}{\rightarrow}

\newcommand{\bs}{\backslash}

\newcommand{\col}{\!:\!}

\newcommand{\fix}{\operatorname{fix}}

\newcommand{\rint}{\operatorname{int}}

\newcommand{\id}{\operatorname{id}}

\newcommand{\is}{\operatorname{Isom}}

\newcommand{\Mob}{\mathop{\rm M\ddot{o}b}\nolimits}

\newenvironment{pf}{\begin{trivlist}\item[]{\bf Proof:\ }}
{\mbox{}\hfill\rule{.08in}{.08in}\end{trivlist}}

\begin{document}

\pagestyle{myheadings}
\markboth{Boris N. Apanasov}{Topological barriers for locally homeomorphic quasiregular mappings} 
\title{Topological barriers for locally homeomorphic quasiregular mappings in 3-space}
\author{Boris N. Apanasov}

\date{}

\maketitle

\begin{abstract}
We construct a new type of locally homeomorphic quasiregular mappings in the 3-sphere and discuss their relation to the M.A.Lavrentiev problem, the Zorich map with an essential singularity at infinity, the Fatou's problem and a quasiregular analogue of domains of holomorphy in complex analysis. The construction of such mappings comes from our construction of non-trivial compact 4-dimensional cobordisms $M$ with symmetric boundary components and whose interiors have complete 4-dimensional real hyperbolic structures. Such locally homeomorphic quasiregular mappings are defined in the 3-sphere $S^3$ as mappings equivariant with the standard conformal action of uniform hyperbolic lattices $\Gamma\subset \is H^3$ in the unit 3-ball and its complement in $S^3$ and with its discrete representation $G=\rho(\Gamma)\subset \is H^4 $. Here $G$ is the fundamental group of our non-trivial hyperbolic 4-cobordism $M=(H^4\cup\Omega(G))/G$  and the kernel of the homomorphism $\rho\col \Ga\rightarrow G$ is a free group $F_3$ on three generators.
\end{abstract}

\maketitle
    \footnotetext[2]{
    \hfil\hfil\hfil}

\section{Introduction}
Liouville's rigidity of spatial conformal geometry shows that conformal mappings in domains in $S^n=\mathbb{R}^n\cup \{\infty\}, n\geq 3$ are restrictions of M\"{o}bius transformations. However this rigidity no longer persists in quasiconformal geometry intensively  studied since 1930s after its introduction by H.Gr\"{o}zsch \cite{Gr} and M.A.Lavrentiev \cite{L1}. First assertions reflecting spatial specifics in this quasiconformal geometry were made by M.A.Lavrentiev \cite{L2}, on removability of some singularities of quasiconformal mappings and on locally homeomorphic mappings in $\mathbb{R}^3$. V.A.Zorich's 1967 solution \cite{Z1} of the last Lavrentiev's problem (the global homeomorphism theorem) shows that locally homeomorphic quasiregular mappings of $\mathbb{R}^n, n\geq 3$, into itself are homeomorphisms of $\mathbb{R}^n$, and thus quasiconformal mappings.

In addition to his famous proof of Lavrentiev's problem Zorich gave an example of a nonsurjective quasiregular mapping $\mathbb{R}^3\rightarrow \mathbb{R}^3$ omitting the origin and having an essential singularity at infinity. This so-called Zorich map is a spatial analogue of the exponential function in $\mathbb{C}$  and is based on P.P.Belinskii's construction of a quasiconformal mapping of a half space $\mathbb{R}^3_{+}$ onto a round solid cylinder. Due to the previous Zorich theorem, the branching of the map (along parallel lines orthogonally intersecting the boundary plane at integer points) cannot be avoided. In a general sense all quasiregular mappings topologically have a branched covering type. Namely by Reshetnyak's theorem, quasiregular mappings are (generalized) branched covers, that is, discrete and open mappings and hence local homeomorphisms modulo an exceptional set of (topological) codimension at least two.
The intensive study of quasiregular mappings, especially after the mentioned Zorich results and conjectures in \cite{Z1}-\cite{Z5}, resulted in a rich theory of quasiregular mappings which is a natural and beautiful generalization of the geometric aspects of the theory of holomorphic functions in the plane to higher dimensions. It is covered by several papers \cite{B}, \cite{BL1}-\cite{BL2}, \cite{DP}, \cite{G1}, \cite{MRV}, \cite{V}, \cite{Z1}-\cite{Z5} and a number of monographs - see \cite{G2}-\cite{G3}, \cite{Re}, \cite{Ri}, \cite{Vu}.

Our motivation stems from three sides of the mentioned Lavrentiev-Zorich assertions: on locally homeomorphic spatial quasiregular (quasimeromorphic) mappings defined in the (almost) whole sphere $S^3=\overline{\mathbb{R}^3}$, surjectivity of such mappings having the whole sphere $S^3$ as the image, and their essential singularities. Despite a relative rigidity of quasiregular mappings without branching, in Theorem \ref{map} we present a new (flexible) way for constructions of such locally homeomorphic quasiregular mappings $S^3\setminus S_*\rightarrow S^3$ defined in the sphere $S^3$ except $S_*$ (a dense subset of the 2-sphere $S^2\subset S^3$) and having the whole sphere $S^3$ as the image. The exceptional subset $S_*$ of the 2-sphere $S^2\subset S^3$ (or a quasi-sphere $S^2_q\subset S^3$) is a countable orbit of a Cantor subset of zero 2-measure. It creates a  barrier (of a topological nature) for continuous extension of our quasiregular mapping since points of $S_*$ are essential singularities of our mapping having no radial limits.

The construction of such quasiregular mappings in $S^3\backslash S_*$ having a dense subset $S_*\subset S^2\subset S^3$ as their barrier
 is heavily based on our construction \cite{A5} of non-trivial compact 4-dimensional cobordisms $M^4$ with symmetric boundary components, which makes it absolutely necessary for understanding to repeat in Section 2 the construction of the corresponding discrete actions and their homomorphisms related to those symmetric 4-cobordisms. The interiors of these 4-cobordisms have complete 4-dimensional real hyperbolic structures and universally covered by the real hyperbolic space $H^4$, while the boundary components of $M^4$ have (symmetric) 3-dimensional conformally flat structures obtained by deformations of the same hyperbolic 3-manifold whose fundamental group $\Gamma$ is a uniform lattice in $\is H^3$. Such conformal deformations of hyperbolic manifolds are well understood after their discovery in \cite{A1}, see \cite{A2}. Nevertheless till recently such "symmetric" hyperbolic 4-cobordisms with described properties were unknown despite our well known constructions of non-trivial hyperbolic homology 4-cobordisms with very assymmetric boundary components - see \cite{AT} and \cite{A2}-\cite{A4}. In \cite{A5} we presented a method of constructing such non-trivial "symmmetric" hyperbolic 4-cobordisms $M^4=H^4/G$ whose fundamental groups
 $\pi_1(M^4)$ act discretely in the hyperbolic 4-space $H^4$ by isometries, $\pi_1(M^4)\cong G\subset \is H^4$, and can be obtained from the hyperbolic 3-lattice $\Gamma\subset \is H^3$ by a homomorphism $\rho\col \Gamma\rightarrow G\subset \is H^4$ with non-trivial kernel (in our construction such kernel of $\rho$ is a free subgroup $F_3\subset\Gamma$ on three generators). In Section 2 we present all necessary details of our construction of such "symmetric" hyperbolic 4-cobordisms and used discrete groups (Theorem \ref{constr} and Proposition \ref{homo}). By using such "symmetric" hyperbolic 4-cobordisms, our locally homeomorphic quasiregular mappings $F$ are defined in the complement $S^3\setminus S^2$ of the 2-sphere $S^2=\{x\in\mathbb{R}^3\col |x|=1\}$ as mappings equivariant with the standard conformal action of uniform hyperbolic lattices $\Gamma\subset \is H^3$ in the unit 3-ball $B^3(0,1)=\{x\in\mathbb{R}^3\col |x|<1\}$ and in the complement in $S^3$ to its closure and with the discrete representation $G=\rho(\Gamma)\subset \is H^4 $. In other words such $\Gamma$-equivariance of our quasiregular mappings $F$ can be described as $F(\Gamma(x))=\rho(\Gamma)(F(x))=G(F(x)), x\in S^3\setminus S^2$. Another essential element of our construction is a direct building in Section 3 of the so called bending quasiconformal homeomorphisms between polyhedra which preserve combinatorial structure of polyhedra and their dihedral angles.

 One may find a resemblance of our construction of locally homeomorphic quasiregular mappings $F$ in Theorem \ref{map} to constructions of O.Martio and U.Srebro \cite{MS1},\cite{MS2} and P.Tukia \cite{T2} of locally homeomorphic quasiregular mappings. However our quasiregular mappings are not automorphic with respect to any discrete M\"{o}bius group as are the mappings in those papers.

 Our construction is also related to the well known open question (Fatou's problem) on the correct analogue for higher-dimensional quasiregular mappings of the Fatou's theorem \cite{F} on radial limits of a bounded analytic function of the unit disc. Though in higher dimensions $n\geq 3$ it is not even known if there exists a bounded $n$-dimensional mapping of the unit ball without any radial limits, cf. \cite{MR}, \cite{V}, there are several results concerning radial limits of mappings of the unit ball. The most recent progress is due to Kai Rajala who in particular proved that radial limits exist for infinitely many points of the unit sphere, see \cite{Ra} and references there for some earlier results in this direction. Considering restriction of our locally homeomorphic quasiregular mapping $F$ in Theorem \ref{map} to the unit ball we see that it is bounded, and its exceptional subset $S_*$ of the boundary unit sphere is a countable orbit of a Cantor subset with Hausdorff dimension $\ln 5/\ln 6 \approx 0.89822444$ (zero 2-measure). All points of $S_*$ are essential singularities of our bounded locally homeomorphic quasiregular mapping having no radial limits.

 \subsection{Acknowledgments}

The author is grateful to Pekka Pankka for fruitful discussions and interest and for organizing with Jang-Mei Wu our informal seminar in Helsinki. We also thank Vladimir A. Zorich for clarifying discussions. We thank the referee for several remarks on the manuscript that led to several improvements.

\section{Non-trivial "symmetric" hyperbolic 4-cobordisms}

Since the construction of the fundamental group
 $\pi_1(M^4)\cong G\subset \is H^4$ of a non-trivial "symmmetric" hyperbolic 4-cobordisms $M^4=H^4/G$  acting discretely in the hyperbolic 4-space $H^4$ is very essential for our construction of a locally homeomorphic quasiregular mapping $F\col S^3\setminus S_*\rightarrow S^3$ having $S_*$ in a quasi-sphere $S^2_q\subset S^3$ as a barrier, we start with a detailed construction of such discrete group $G\subset \is H^4$ and the corresponding discrete representation $\rho\col \Gamma \rightarrow G$ of a uniform hyperbolic lattice $\Gamma\subset \is H^3$ from our paper \cite{A5}.

These discrete groups $G$ and $\Gamma$ negatively answer a conjecture: \textit{If one had a hyperbolic 4-cobordism $M^4$ whose boundary components $N_1$ and $N_2$ are highly (topologically and geometrically) symmetric to each other it would be in fact an h-cobordism, possibly not trivial, i.e. not homeomorphic to the product of $N_1$ and  the segment $[0,1]$.}

Namely the boundary components $N_1$ and $N_2$ of $M^4=M(G)=\{H^4\cup \Omega(G)\}/G$ are covered by the discontinuity set $\Omega (G)\subset S^3$ of $G$ with two connected components $\Omega_1$ and $\Omega_2$, where the conformal action of $G=\rho(\Gamma)$ is symmetric and has contractible fundamental polyhedra $P_1$ and $P_2$ of the same combinatorial type allowing to realize them as a compact polyhedron $P_0$ in the hyperbolic 3-space, i.e.
the dihedral angle data of these polyhedra satisfy the Andreev's conditions \cite{An1}. Nevertheless this geometric symmetry of boundary components of our hyperbolic 4-cobordism $M(G)$) is not enough to ensure
that the group $G=\pi_1(M^4)$ is quasi-Fuchsian and our 4-cobordism $M$ is trivial.

 Here a Fuchsian
group $\Ga\subset\is H^3\subset \is H^4$ conformally acts in the 3-sphere $S^3=\p H^4$ and preserves a round ball $B^3\subset S^3$ where it acts as a cocompact discrete group of isometries of $H^3$. Due to the Sullivan
structural stability (see Sullivan \cite{S} for $n=2$ and Apanasov \cite{A2}, Theorem 7.2), the
space of quasi-Fuchsian representations of a hyperbolic lattice $\Ga\subset\is H^3$ into $\is H^4$
is an open connected component of the Teichm\"uller space of $H^3/\Ga$ or the variety of conjugacy classes of discrete representations
$\rho\col \Ga\ra\is H^4$. Points in this (quasi-Fuchsian) component correspond to trivial hyperbolic 4-cobordisms
$M(G)$ where the discontinuity set $\Om(G)=\Om_1\cup\Om_2\subset S^3=\p H^4$ is the union of two topological 3-balls $\Om_i$, $i=1,2$, and
$M(G)$ is homeomorphic to the product of $N_1$ and the closed interval $[0,1]$.

To simplify the situation we may consider the hyperbolic 4-cobordisms $M(\rho(\Ga))$ corresponding to uniform hyperbolic lattices $\Ga\subset \is H^3$ generated by reflections (or cobordisms related to their finite index subgroups). Natural inclusions of these lattices into $\is H^4$ act at infinity $\p H^4 = S^3$ as Fuchsian groups
$\Ga\subset \Mob (3)$ preserving a round ball in the 3-sphere $S^3$. In this case the  above conjecture can be reformulated as the
following question on the M\"obius action of corresponding reflection groups $G=\rho(\Ga)\subset \is H^4$ on the 3-sphere $S^3=\p H^4$:

\begin{quest}\label{quest} Is any discrete M\"obius group $G$ generated by finitely many
reflections with respect to spheres $S^2\subset S^3$ and whose
fundamental polyhedron $P(G)\subset S^3$ is the union of two contractible polyhedra
$P_1, P_2\subset S^3$ of the same combinatorial type (with equal corresponding dihedral angles) quasiconformally conjugate
in the sphere $S^3$ to some Fuchsian group preserving a round ball $B^3\subset S^3$?
\end{quest}

 Our construction of the mentioned discrete groups $\Gamma$ and $G=\rho(\Gamma)$ gives a negative answer to this question and proves the following (see Apanasov\cite{A5}):

\begin{theorem}\label{constr}
There exists a discrete M\"obius group $G\subset \Mob (3)$ on the 3-sphere $S^3$ generated by finitely many reflections such that:
\begin{enumerate}
\item Its discontinuity set $\Om(G)$ is the union
of two invariant components $\Om_1$, $\Om_2$;
\item Its fundamental polyhedron $P\subset S^3$ has two contractible components $P_i\subset\Om_i$, $i=1,2$,
having the same combinatorial type (of a compact hyperbolic polyhedron $P_0\subset H^3$);
\item For the uniform hyperbolic lattice $\Ga\subset\is H^3$ generated by reflections in sides of the hyperbolic
polyhedron $P_0\subset H^3$ and acting on the sphere $S^3=\p H^4$ as a discrete Fuchsian group  $i(\Ga)\subset \is H^4=\Mob(3)$ preserving a round ball $B^3$ (where $i\col\is H^3\subset\is H^4$ is the natural inclusion),
the group $G$ is its image under a homomorphism $\rho\col\Ga\ra G$  but it is not quasiconformally (topologically) conjugate in $S^3$ to $i(\Ga)$.
\end{enumerate}
\end{theorem}

\begin{pf}
For our construction of the desired M\"obius group $G\subset \Mob (3)$ generated by reflections it is
enough to define its finite collection $\Sa$ of reflecting 2-spheres
$S_i\subset S^3$, $1\leq i \leq N$. As the first four spheres  we consider mutually orthogonal spheres
centered at the vertices of a regular tetrahedron in
${\mathbb R}^3$. Let $B=\bigcup_{1\leq i\leq 4} B_i$ be the union of the closed balls bounded by these four spheres,
and let $\p B$ be its boundary (a topological 2-sphere) having four vertices which are the intersection points of four triples of our spheres.
Applying a  M\"obius transformation in $S^3\cong{\mathbb R}^3\cup\{\infty\}$, we may
assume that the first three spheres $S_1, S_2$ and  $S_3$
correspond to the coordinate planes $\{x\in {\mathbb R}^3\col x_i=0\}$, and
 $S_4=S^2(0,R)$ is the round sphere of some radius $R>0$ centered at the origin. The value of the radius $R$ will be determined later.

On the topological 2-sphere $\p B$ with four vertices we consider a simple closed loop $\alpha\subset \p B$
which does not contain any of our vertices and which symmetrically separates two pairs of these vertices from each other
 as the white loop does on the tennis ball shown in Figure \ref{fig6}. This white loop
$\alpha$ can be considered as the boundary of a topological 2-disc $\sigma$ embedded in the complement
   $D=S^3\setminus B$ of our four balls. Our geometric construction needs a detailed
   description of such a 2-disc $\sa$ and its boundary loop $\al=\p \sa$ obtained as it is shown in
   Figure \ref{fig7}.

The desired disc $\sa\subset D=S^3\setminus B$ can be described as the boundary in the domain $D$ of the union
of a finite chain of adjacent blocks $Q_i$ (regular cubes) with disjoint interiors whose centers
lie on the coordinate planes $S_1$ and $S_2$ and whose sides are parallel to the coordinate planes.
This chain starts from the unit cube whose center lies
in the second coordinate axis, in $e_2\cdot \mathbb R_{+}\subset S_1\cap S_3$. Then our chain goes up through small adjacent cubes centered in the coordinate plane $S_1$, at some point changes its direction to the horizontal one toward the third coordinate axis, where it turns its horizontal direction by a right angle again (along the coordinate plane $S_2$), goes toward the vertical line passing through the second unit cube centered in
$e_1\cdot \mathbb R_{+}\subset S_2\cap S_3$, then goes down along that vertical line and finally ends
at that second unit cube, see Figure \ref{fig7}. We will define the size of small cubes $Q_i$ in our block chain
and the distance of the centers of two unit cubes to the origin in the next step of our construction.
\begin{figure}
\centering
\epsfxsize=8cm
\epsfbox{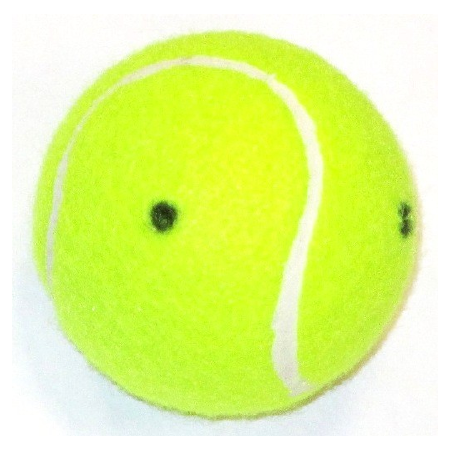}
\caption{White loop separating two pairs of vertices on a tennis ball.}
\label{fig6}
\end{figure}

\begin{figure}
\centering
\epsfxsize=9cm
\epsfbox{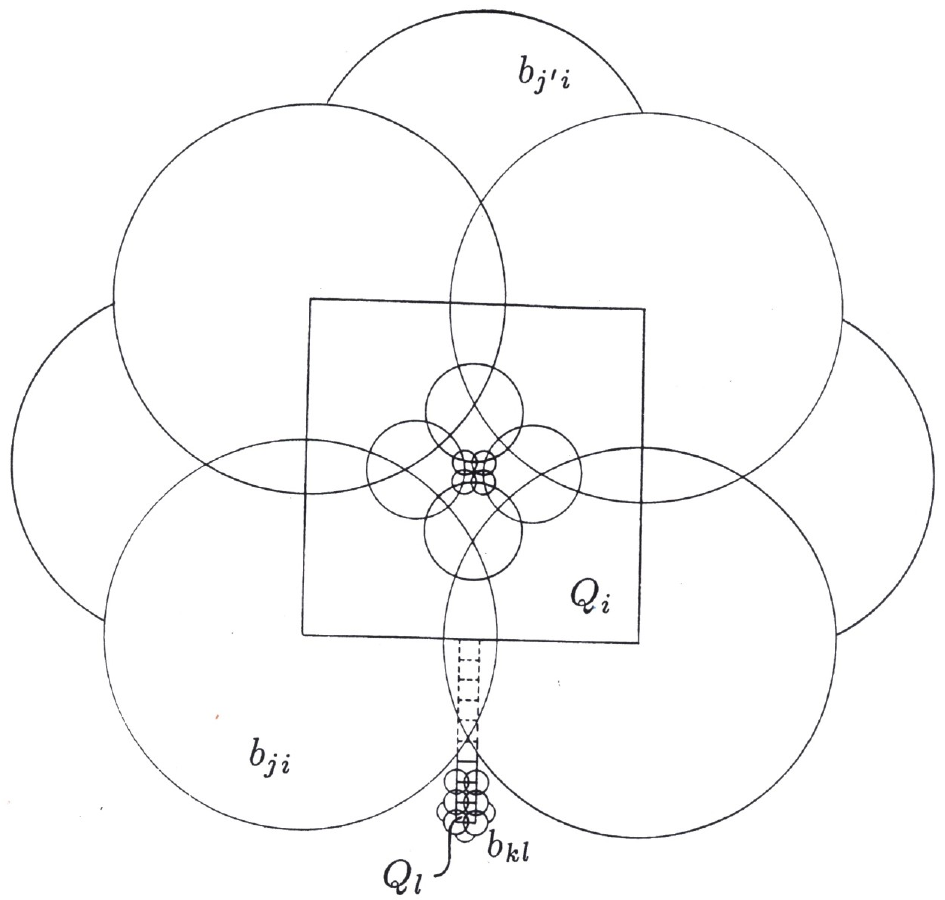}
\caption{Big and small cube sizes and ball covering}
\label{fig4}
\end{figure}

\begin{figure}
\centering
\epsfxsize=14cm
\epsfbox{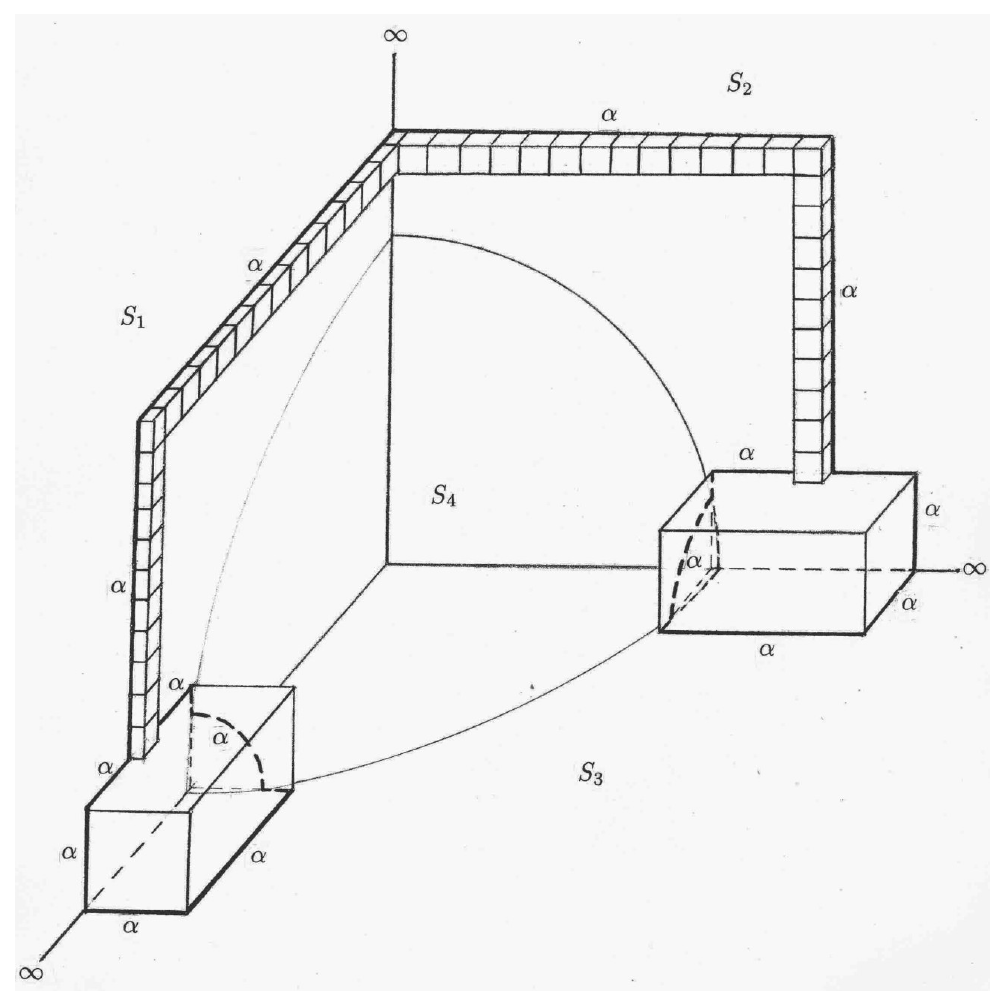}
\caption{Configuration of blocks and the white loop $\alpha\subset \p B$.}
\label{fig7}
\end{figure}
Let us consider one of our cubes $Q_i$, i.e. a block of our chain, and let $f$ be its square side having a nontrivial intersection with our 2-disc
$\sa\subset D$.
For that side $f$ we consider spheres $S_j$ centered at its vertices and
having a radius such that each two spheres centered at the ends of an edge of $f$ intersect each other with angle $\pi/3$.
In particular, for the unit cubes such spheres have radius $ \sqrt{3}/3 $. From such defined spheres we select those spheres that have centers
in our domain $D$ and then include them in the collection $\Sa$ of reflecting spheres.
Now we define the distance of the centers of our big (unit) cubes to the origin. It is determined by the condition that the sphere $S_4=S^2(0,R)$
is orthogonal to the sphere $S_j\in\Sa$ centered at the vertex of such a cube closest to the origin.

As in Figure \ref{fig4}, let $f$ be a square side of one of our cubic blocks  $Q_i$ having a nontrivial intersection $f_{\sa}=f\cap\sa$ with our
2-disc $\sa\subset D$.
We consider a ring of four spheres $S_i$ whose centers are interior points of $f$ which lie outside of the four
previously defined spheres $S_j$ centered at vertices of $f$ and such that each sphere $S_i$ intersects two adjacent spheres $S_{i-1}$ and $S_{i+1}$
(we numerate spheres $S_i$ mod 4) with angle $\pi/3$. In addition these spheres $S_i$
are orthogonal to the  previously defined ring of bigger spheres $S_j$, see Figure \ref{fig4}.
From such defined spheres $S_i$ we select those spheres that have nontrivial
intersections with our domain $D$ outside the previously defined spheres $S_j$, and then include them in the collection $\Sa$ of reflecting spheres.
If our side $f$ is not the top side of one of the two unit cubes we add another sphere $S_k\in\Sa$.
It is centered at the center of this side $f$ and is orthogonal to the four previously defined spheres $S_i$ with centers in $f$, see Figure \ref{fig4}.

Now let $f$ be the top side of one of the two unit cubes of our chain. Then, as before,
we consider another ring of four spheres $S_k$. Their centers are interior points of $f$, lie outside of the four
previously defined spheres $S_i$ closer to the center of $f$ and such that each sphere $S_k$ intersects two adjacent spheres $S_{k-1}$ and $S_{k+1}$
(we numerate spheres $S_k$ mod 4) with angle $\pi/3$. In addition these new four spheres $S_k$
are orthogonal to the previously defined ring of bigger spheres $S_i$, see Figure \ref{fig4}.
 We note that the centers of these four new spheres $S_k$ are vertices of a small square
$f_s\subset f$ whose edges are parallel to the edges of $f$, see Figure \ref{fig4}. We set this square $f_s$ as the bottom side of the small cubic box adjacent to the unit one.
This finishes our definition of the family of twelve round spheres whose interiors cover the square ring $f\bs f_s$ on the top side of one of the two unit cubes
in our cube chain and tells us which two spheres among the four new defined spheres $S_k$ were already included
in the collection $\Sa$ of reflecting spheres
(as the spheres $S_j\in\Sa$ associated to small cubes in the first step).

This also defines the size of small cubes in our block chain. Now we can vary the remaining free parameter $R$
(which is the radius of the sphere $S_4\in\Sa$) in order to make two horizontal rows of small blocks with centers in $S_1$ and $S_2$, correspondingly,
to share a common cubic block centered at a point in $e_3\cdot \mathbb R_{+}\subset S_1\cap S_2$, see Figure \ref{fig7}.

The constructed collection $\Sa$ of reflecting spheres $S_j$ bounding round balls $B_j$, $1\leq j\leq N$, has the following properties:

\begin{enumerate}
\item The closure of our 2-disc $\sa\subset D$ is covered by balls $B_j$: $\bar{\sigma}\subset \rint\bigcup\limits_{j\geq5}^{N}B_j$;
\item Any two spheres $S_j, S_{j'}\in\Sa$ 
either are disjoint or intersect with angle $\pi/2$ or  $\pi/3$;
\item The complement of all balls, $S^3\setminus\bigcup\limits_{j=1}^{N}B_j$ is the union of two contractible
polyhedra $P_1$ and $P_2$ of the same combinatorial type.
\end{enumerate}

Therefore we can use the constructed collection $\Sa$ of reflecting spheres $S_i$ to define a discrete group $G=G_{\Sa}\subset \Mob(3)$
generated by $N$ reflections in spheres $S_j\in \Sa$. The fundamental polyhedron $P=P_1\cup P_2\subset S^3$ for the action of this discrete reflection group
$G$ on the sphere $S^3$ is the union of two connected polyhedra
$P_1$ and $P_2$ which are disjoint topological balls. So the discontinuity set $\Omega(G)\subset S^3$ of $G$ consists
of two invariant connected components $\Om_1$ and $\Om_2$:

\begin{equation}\label{comp}
\Omega(G)=\bigcup_{g\in G} g(\bar{P})=\Om_1\cup\Om_2\,,\quad \Om_i=\bigcup_{g\in G} g(\bar{P_i})\,, \quad i=1,2.
\end{equation}

\begin{lemma}\label{h-body}
The splitting of the discontinuity set $\Om\subset S^3$ of our discrete reflection group $G=G_{\Sa}\subset \Mob(3)$ into $G$-invariant components\, $\Om_1$ and $\Om_2$ in (\ref{comp}) defines a Heegaard splitting of the 3-sphere $S^3$ of infinite genus with ergodic word hyperbolic group $G$ action on the separating boundary $\La(G)$ which is quasi-self-similar in the sense of Sullivan.
\end{lemma}

\begin{pf}
In fact, despite the contractibility of polyhedra $P_1$ and $P_2$ both components $\Om_1$ and $\Om_2$ are not simply connected and even are mutually linked.
To show this it is enough to see that the union
of the bounded polyhedron $\bar{P_1}$ (inside of our block chain) and its image $g_3(\bar{P_1})$ under the reflection $g_3$ with respect to the plane
$S_3$ has
a non-contractible loop $\be_1$ which represents a non-trivial element of the fundamental group $\pi_1(\Om_1)$.
This loop is linked with the loop $\be_2$ in the unbounded component $\Om_2$ which goes around
$\bar{P_1}\cup g_3(\bar{P_1})$ and represents a non-trivial element of the fundamental group $\pi_1(\Om_2)$.

\begin{figure}
\centering
\includegraphics[width=13cm]{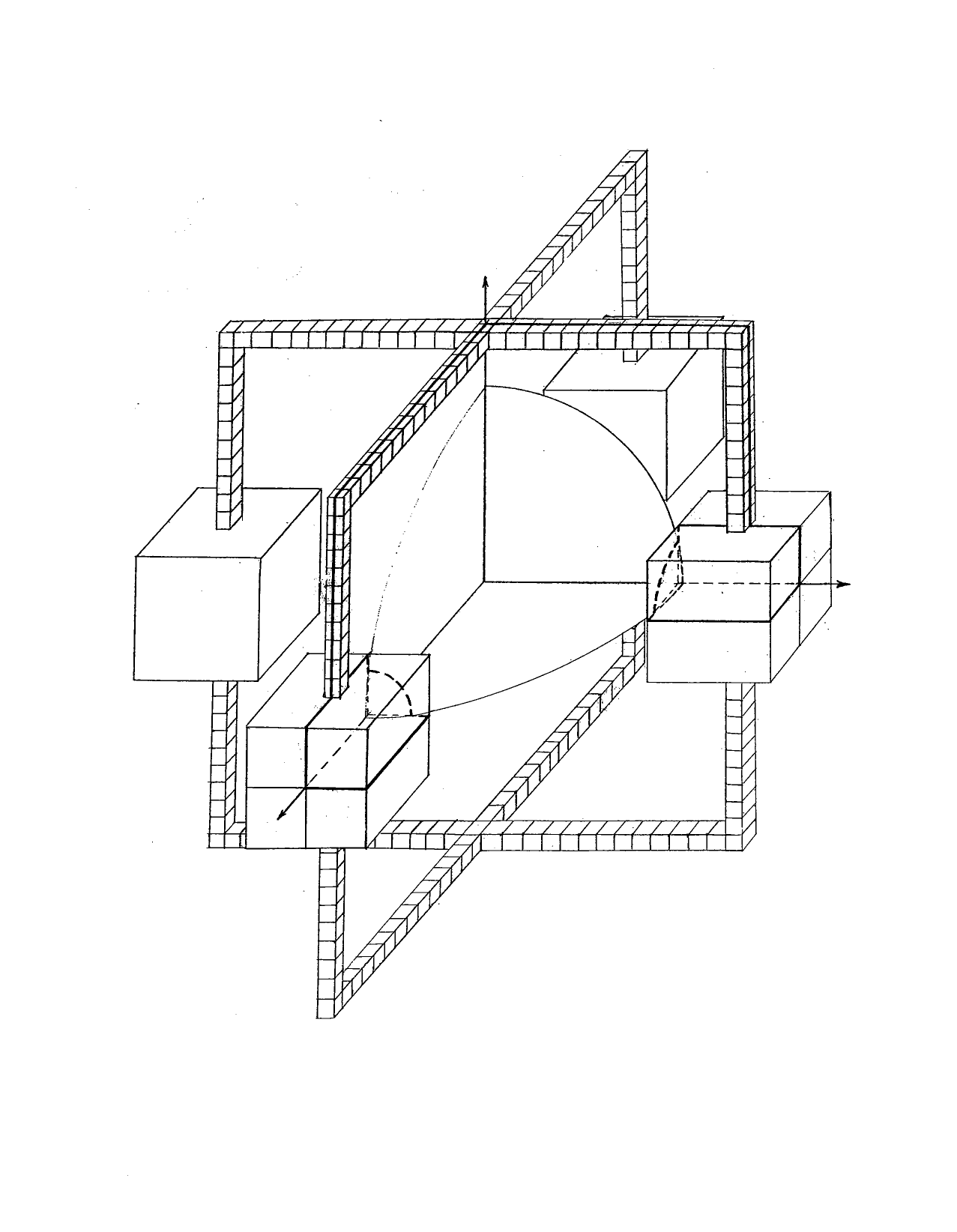}
\caption{Handlebody obtained by the first 3 reflections of the cube chain.}
\label{fig8}
\end{figure}
This fact is illustrated by Figure \ref{fig8} where one can
see a handlebody obtained from our initial chain of building blocks in Figure \ref{fig7} by the union of the images of this block chain by first generating reflections in the group $G$ (in $S_1, S_2$ and $S_3$). Then our non-contractible loop $\be_1\subset\Om_1$ lies inside of this handlebody in Figure \ref{fig8}  and is linked with the second loop $\be_2\subset\Om_2$ which goes around one of the handles of the handlebody in Figure \ref{fig8}. The resulting
handlebodies $\Omega_1$ and  $\Omega_2$ are the unions of the corresponding images $g(\bar{P_i})$ of the polyhedra $\bar{P_1}$ and $\bar{P_2}$, so they have infinitely many mutually linked handles. Their fundamental groups $\pi_1(\Om_1)$ and  $\pi_1(\Om_2)$ have infinitely many generators, and some of those generators correspond to the group $G$-images of the linked loops
$\be_1\subset\Om_1$ and $\be_2\subset\Om_2$. The limit set
$\La(G)$ is the common boundary of $\Om_1$ and $\Om_2$. Since the group $G\subset\Mob(3)$ acts on the hyperbolic 4-space $H^4$, $\p H^4=S^3$,
as a convex cocompact isometry group, its action on the limit set $\La(G)$ is ergodic. Moreover, the common boundary
$\La(G)$ of the handlebodies $\Omega_1$ and  $\Omega_2$ is quasi-self-similar in the sense of Sullivan, that is
each arbitrary small piece of $\La(G)$ can be expanded to a standard size and then mapped into $\La(G)$
by a $K$-quasi-isometry. More precisely, there are uniform constants $K$ and $r_0$ such that, for any $x\in \La(G)$
and  for any ball $B(x,r)$ centered at $x$ with radius $r$, $0<r<r_0$, there exists a
$K$-quasi-isometric bijection $f$,
 \begin{eqnarray}\label{qi}
 f\,:\,\, \frac {\La(G)\cap B(x,r)}{r} \hookrightarrow \La(G)
\end{eqnarray}
which distorts distances in the interval between $1/K$ and $K$. In other words,
the distortion of an unlimited {\it ``microscoping"} (\ref{qi}) of the limit set $\La(G)$ can be uniformly bounded, see Corollary 2.66 in Apanasov \cite{A2}.
\end{pf}

To finish the proof of Theorem \ref{constr} we notice that the combinatorial type (with magnitudes of dihedral angles) of the bounded component $P_1$ of the fundamental polyhedron $P\subset S^3$ coincides with
the combinatorial type of its unbounded component $P_2$. Applying Andreev's theorem on 3-dimensional hyperbolic polyhedra \cite{An1},
one can see that there exists a compact hyperbolic polyhedron $P_0\subset H^3$ of the same combinatorial type with the same dihedral angles ($\pi/2$ or $\pi/3$).
So one can consider a uniform hyperbolic lattice $\Ga\subset\is H^3$ generated by reflections in sides of the hyperbolic
polyhedron $P_0$. This hyperbolic lattice $\Ga$ acts in the sphere $S^3$ as a discrete co-compact Fuchsian group $i(\Ga)\subset \is H^4=\Mob(3)$
(i.e. as the group
$i(\Ga)\subset \is H^4$ where $i\col\is H^3\subset\is H^4$ is the natural inclusion)
preserving a round ball $B^3$ and having its boundary sphere $S^2=\p B^3$ as the limit set. Obviously there is no self-homeomorphism of the sphere $S^3$
conjugating the action of the groups $G$ and $i(\Ga)$ because the limit set $\La(G)$ is not a topological 2-sphere. So the constructed group $G$ is not a
quasi-Fuchsian group.

One can construct a natural homomorphism $\rho\col \Ga\ra G$, $\rho\in\scr_3(\Ga)$, between these two Gromov hyperbolic groups $G\subset \is H^4$ and $\Ga\subset \is H^3$
defined by the correspondence between sides of
the hyperbolic polyhedron $P_0\subset H^3$ and reflecting spheres $S_i$ in the collection $\Sa$ bounding the fundamental polyhedra $P_1$ and $P_2$.

\end{pf}

\begin{prop}\label{homo}
The homomorphism $\rho\in\scr_3(\Ga)$, $\rho\col \Ga\ra G$, in Theorem \ref{constr} is not an isomorphism. Its kernel $\ker(\rho)=\rho^{-1}(e_G)$ is a
free rank 3 subgroup $F_3\lhd\Ga$.
\end{prop}

\begin{pf}
The homomorphism $\rho$ cannot be an isomorphism since its kernel $\rho^{-1}(e_G)$ is not trivial,
$\rho^{-1}(e_G)\neq \{e_{\Ga}\}$. In fact this kernel is a
free rank 3 group $F_3=\langle x, y, z\rangle$ generated by three hyperbolic translations $x, y, z \in\Ga$. The first hyperbolic translation  $x=a_1b_1$ in $H^3$ is the composition of reflections $a_1$ and $b_1$ in two disjoint hyperbolic planes $H_1, H'_1\subset H^3$ containing those two 2-dimensional faces of
the hyperbolic polyhedron $P_0$ that correspond to two sides of the polyhedron $P_1$ which are disjoint parts of the sphere $S_4$.  The second
 hyperbolic translation  $y=a_2b_2$ in $H^3$ is the composition of reflections $a_2$ and $b_2$ in two disjoint hyperbolic planes $H_2, H'_2\subset H^3$  containing those two 2-dimensional faces of the hyperbolic polyhedron $P_0$ that correspond to two sides of the polyhedron $P_1$ which are disjoint parts of the sphere $S_3$.
 And the third generator $z$ is a hyperbolic translation in $H^3$ which is $a_1$-conjugate of $y$, $z=a_1ya_1$. The fact that these hyperbolic 2-planes $H_1$ and $H'_1$ (correspondingly, the 2-planes $H_2$ and $H'_2$) are disjoint follows from Andreev's result \cite{An2} on sharp angled hyperbolic polyhedra. Restricting our homomorphism $\rho$ to the subgroup of $\Ga$ generated by reflections $a_1, a_2, b_1, b_2\in \Ga$, we can formulate its properties as the following statement in combinatorial group theory:

 \begin{lemma}\label{ker}
Let $A=\langle a_1, a_2 \mid a_1^2,\, a_2^2,\, (a_1a_2)^2\rangle $ $\cong$ $B=\langle b_1, b_2 \mid b_1^2,\, b_2^2,\, (b_1b_2)^2\rangle $ $\cong$
$C=\langle c_1, c_2 \mid c_1^2,\, c_2^2,\, (c_1c_2)^2\rangle \cong \mathbb{Z}_2 \times \mathbb{Z}_2$,
and let $\varphi\col A\ast B\ra C$ be a homomorphism of the free product $A\ast B$ into $C$ such that $\varphi(a_1)=\varphi(b_1)=c_1$ and
 $\varphi(a_2)=\varphi(b_2)=c_2$. Then the kernel $\ker(\varphi)=\varphi^{-1}(e_C)$ of $\varphi$ is a free rank 3 subgroup $F_3\lhd A\ast B$  generated by elements
 $x=a_1b_1$, $y=a_2b_2$ and $z=a_1a_2b_2a_1=a_1ya_1$.
\end{lemma}

\begin{pf}
 It is obvious that $K_0=\langle x, y, z \rangle$ is a subgroup
in $\ker(\varphi)$.  From the definition of $\varphi$  on the generators of $A\ast B$ it is also clear
that $\ker(\varphi)$ is $\langle\langle x, y \rangle\rangle$, the normal closure of elements
$x=a_1b_1$ and $y=a_2b_2$. Therefore in order to prove that $K_0=\ker(\varphi)$ it is enough to show
that $K_0$ contains all elements which are conjugate in $A\ast B$ to $x$ and $y$.

As any element $w \in A\ast B$ is a product of generators of $A\ast B$,
the conjugation by any such $w$ may be regarded as a consequent conjugation
by the generators of $A \ast B$. So, it is enough to prove
that $K_0$ contains any element conjugate to $x, y, z$ by $a_1, a_2, b_1, b_2$.

In fact it is easy to verify that

$$\begin{array}{lll}

a_1^{-1}xa_1=x^{-1},& \,a_1^{-1}ya_1=z, & \,a_1^{-1}za_1=y,\\

a_2^{-1}xa_2=zxy^{-1},& \,a_2^{-1}ya_2=y^{-1}, & \,a_2^{-1}za_2=z^{-1},\\

b_1^{-1}xb_1=x^{-1}, & \,b_1^{-1}yb_1=x^{-1}zx, & \,b_1^{-1}zb_1=x^{-1}yx,\\

b_2^{-1}xb_2=y^{-1}zx, & \,b_2^{-1}yb_2=y^{-1}, & \,b_2^{-1}zb_2=y^{-1}z^{-1}y.
\end{array}$$

Now we should show that the elements $x, y$ and $z$ form a free basis for $K_0$.
Let us check that any reduced word $w(x,y,z)$ represents a nontrivial element
of $A\ast B$.  By a ``letter'' we mean any of symbols $x^{\pm 1}, y^{\pm 1}, z^{\pm 1}$.
We claim that for the element $g$ represented by a reduced word
$w(x,y,z)$ the following holds: the last syllable of $g$ written in the normal form
is always equal to the last syllable of the last letter of the word $w$ (and so is nontrivial)
except for the case when the last letter is $x$ and the preceding letter is $z$ (cf. \cite{RS}, \S 4.1). In this case the last syllable
equals $b_1b_2$ (and so it is also nontrivial).
Besides that, if the last letter of $w$ is $z$
then the two last syllables of $w$ are the ones of $z$.

This statement can be easily verified by induction on the length of $w(x,y,z)$. So it obviously implies nontriviality of
the element $g$ represented by $w(x,y,z)$.
\end{pf}

\begin{figure}
\centering
\epsfxsize=6cm
\epsfbox{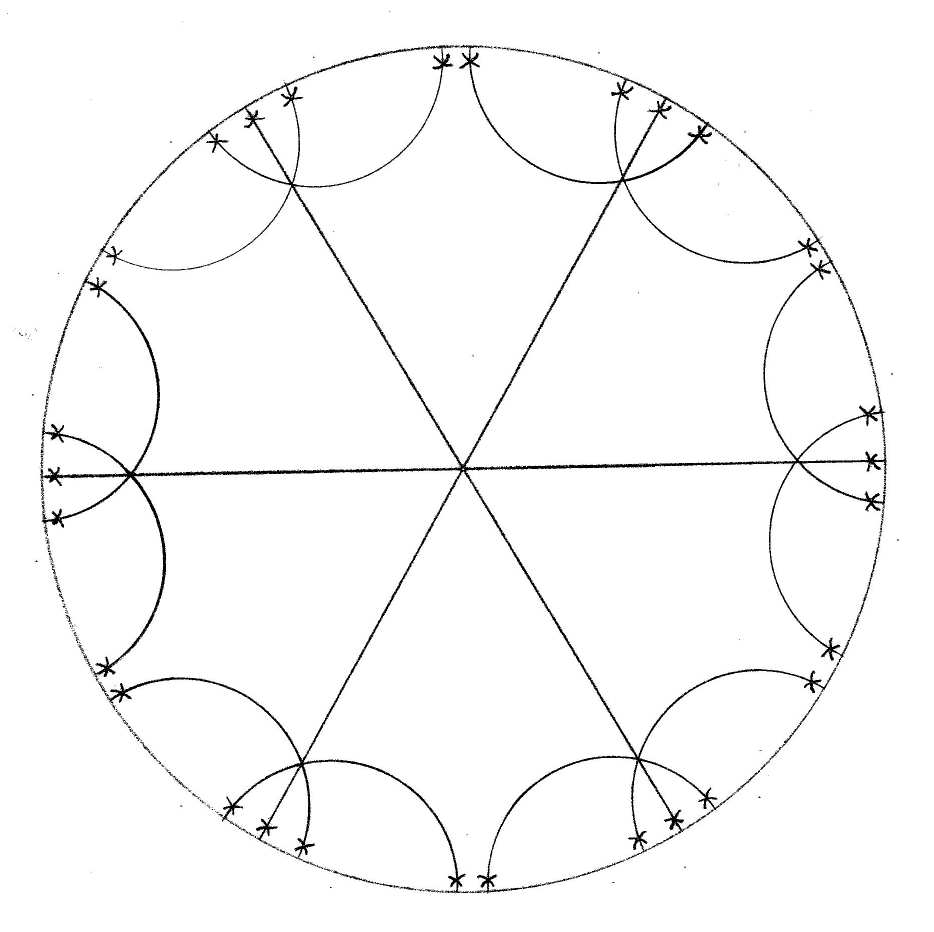}
\caption{Cayley graph of the free group $F_3$.}
\label{fig10}
\end{figure}
Now the claim that $\ker(\rho)\subset\Ga$ is a free rank 3 subgroup
$F_3=\langle x, y, z\rangle$ generated by our hyperbolic translations $x, y, z\in \is H^3$ follows directly from Lemma \ref{ker}, which completes the proof of Proposition \ref{homo}.
\end{pf}

Therefore the configuration of reflecting spheres $S_j\subset \Sa$ shows that one can deform our discrete co-compact Fuchsian group
$i(\Ga)\subset \is H^4=\Mob(3)$ preserving a round ball $B^3\subset S^3$ into the group $G\subset \is H^4$ by continuously moving two pairs
of reflecting 2-spheres of the Fuchsian group $i(\Ga)$ corresponding to the pairs of hyperbolic planes $H_1, H'_1\subset H^3$ and $H_2, H'_2\subset H^3$
into the reflecting spheres $S_4$ and $S_3$ while keeping all dihedral angles unchanged.

\begin{remark}
A simple but important observation is that in our construction we can change the unit round ball $B(0,1)\subset S^3$ to any quasiball $B_q\subset S^3$ bounded by a quasisphere $S_q$ where our hyperbolic lattice $\Gamma\subset \is H^3$ acts conformally and cocompactly. This action is a quasi-Fuchsian group action corresponding to a point in the space of quasi-Fuchsian representations of $\Ga\subset \is H^3$ into $\is H^4$ (an open connected component of the Teichm\"uller space of conformally flat structures on $H^3/\Ga$), cf. \cite{A1, A2}.
\end{remark}

\section{Bending homeomorphisms between polyhedra}

In this section we construct quasiconformal homeomorphisms $\phi_1\col P_1\rightarrow P_0$ and $\phi_2\col P_2\rightarrow \widehat{P_0}$ between components $P_i, i=1,2$, of the fundamental polyhedron $P\subset\Omega(G)\subset S^3$ for the group $G$ and the corresponding components $P_0$ and $\widehat{P_0}$ of the fundamental polyhedron for conformal action in $S^3$ of our hyperbolic lattice $\Gamma\subset\is H^3$ from Theorem \ref{constr}. These mappings $\phi_i$ are compositions of finitely many elementary "bending homeomorphisms", map faces to faces, and preserve the combinatorial structure of the polyhedra and their corresponding dihedral angles.

First we observe that to each cube $Q_j, 1\leq j\leq m$, used in the previous section for our construction of the group $G$ (see Figure \ref{fig7} and Figure \ref{fig8}), we may associate a round
ball $B_j$ centered at the center of the cube $Q_j$ and such that its boundary sphere is orthogonal to the reflection spheres $S_i$ from our generating family $\Sigma$ whose centers are at vertices of the cube $Q_j$. In particular for the unit cubes $Q_1$ and $Q_m$, the reflection spheres $S_i$ centered at their vertices have radius $\sqrt{3}/3$, so the balls $B_1$ and $B_m$ (whose boundary spheres are orthogonal to those corresponding reflection spheres $S_i$) should have radius $\sqrt{5/12}$. Also we add another extra ball $B^3(0,R)$ (which we consider as two balls $B_0$ and $B_{m+1}$) whose boundary is the reflection sphere
$S^2(0,R)=S_4\in\Sigma$ centered at the origin and orthogonal to the closest reflection spheres $S_i$ centered at vertices of two unit cubes $Q_1$ and $Q_m$. Our different enumeration of this ball will be used when we consider different faces of our fundamental polyhedron $P_1$ lying on that reflection sphere $S_4$.

Now for each cube $Q_j, 1\leq j\leq m$, we may associate a discrete subgroup $G_j\subset G\subset\Mob(3)\cong\is H^4$ generated by reflections in the spheres  $S_i\in\Sigma$ associated to that cube $Q_j$ - see our construction in Theorem \ref{constr}. One may think about such a group $G_j$ as a result of quasiconformal bending deformations (see \cite{A2}, Chapter 5) of a discrete M\"{o}bius group preserving the round ball $B_j$ associated to the cube $Q_j$ (whose center coincides with the center of the cube $Q_j$). As the first step in such deformations, let  us define two
quasiconformal ``bending'' self-homeomorphisms of $S^3$, $f_1$ and
$f_{m+1}$, preserving the balls $B_1,\ldots,B_m$ and the set
of their reflection spheres $S_i$, $i\neq 4$, and transferring $\p B_0$ and $\p B{m+1}$ into 2-spheres
orthogonally intersecting $\p B_1$ and $\p B_m$ along round circles $b_1$ and $b_{m+1}$,
respectively.

\begin{figure}
\centering
\epsfxsize=10cm
\epsfbox{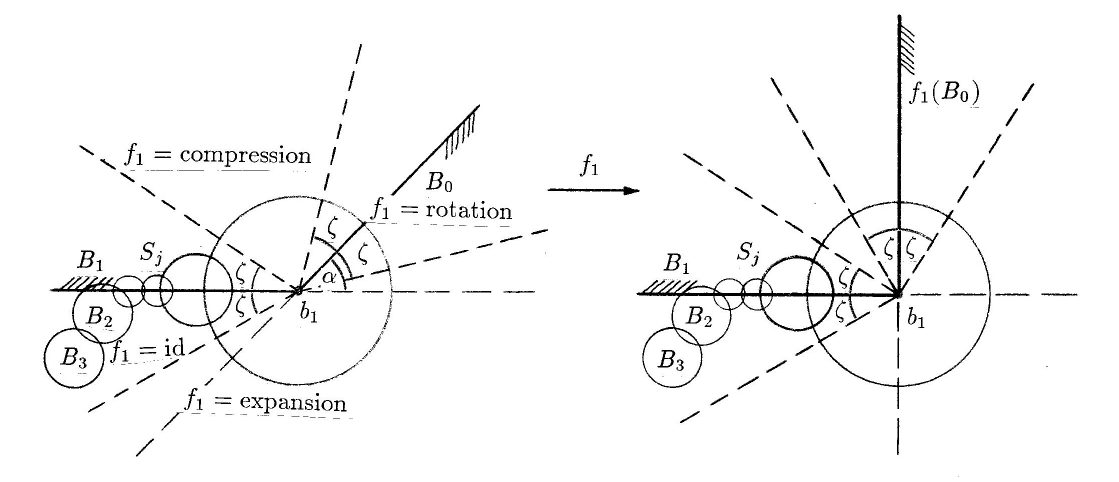}
\caption{Elementary bending homeomorphism $f_1$.}
\label{fig11}
\end{figure}

To construct the bending  $f_1$ ($f_{m+1}$ is similar), we
may assume that the balls $B_0$ and $B_1$ are half-spaces with boundary planes $\p B_0$
and $\p B_1$ and such that
$b_1=\{x\in\mathbb{R}^3\col x_1=x_2=0\}$ is their intersection line.  From our construction of the group
$G$, we have that the dihedral angle of the intersection $B_0\cap B_1$
has a magnitude $\alpha$ , $0<\alpha<\pi/2$, and there exists a dihedral
angle $V_1\subset\mathbb{R}^3$ with the edge $b_1$ and magnitude $2\zeta$,
where
$0<\zeta<\pi/4$ and $\alpha<\pi-2\zeta$, such that $V_1$ contains all the reflection spheres
in $\Sigma$ disjoint from $b_1$.  Let us assume the natural complex structure in
the orthogonal to $b_1$ plane $\mathbb{R}^2=\{x\in\mathbb{R}^3\col x_3=0\}$.  Then the
quasiconformal homeomorphism $f_1\col S^3\ra S^3$ is described by its
restriction to this plane $\mathbb{C}=\mathbb{R}^2$ (where $-\pi<\arg z\leq \pi$ is the
principal value of the argument of $z\in\mathbb{C}$) as follows, see Figure \ref{fig11}:

\begin{eqnarray}\label{bend-1}
f_1(z)=
\begin{cases}
z & \text{if}\ |\arg z|\geq\pi-\zeta\\
z\exp\left(i(\frac{\pi}{2}-\alpha)\right) & \text{if}\ |\alpha-\arg z|\leq\zeta\\
z\exp\left(i(\frac{\pi}{2}-\alpha)(1-\frac{\arg z-\zeta}{\pi-2\zeta})\right) & \text{if}\ \alpha+\zeta<\arg z<\pi-\zeta\\
z\exp\left(i(\frac{\pi}{2}-\alpha)(1+\frac{\zeta+\arg z}{\pi-2\zeta})\right) & \text{if}\ \zeta-\alpha<\arg z<\alpha-\zeta
\end{cases}
\end{eqnarray}

\noindent
We remark that $f_1=\id$ in $V_1$ and hence it is the identity on all reflection spheres
$S_i\in\Sigma$ disjoint from $b_1=\p B_0\cap\p B_1$. Also all spheres $S_k\in\Sigma$
intersecting $b_1$ and the exterior dihedral angles of their
intersections  with other
spheres $S_i$ are still invariant with respect to $f_1$.

 In the next steps in our bending deformations, for two adjacent cubes  $Q_{j-1}$ and  $Q_j$, let  us denote $G_{j-1,j}\subset G$ the subgroup  generated by reflections with respect to the spheres $S_i\subset\Sigma$ centered at common vertices of these cubes. This subgroup preserves the round circle $b_j=b_{j-1,j}=\p B_{j-1}\cap\p B_j$. This shows that our group $G$ is a result of the so called "block-building construction" (see \cite{A2}, Section 5.4) from the block groups $G_j$ by sequential amalgamated products:
\begin{eqnarray}\label{amalgam}
G=G_1\underset{G_{1,2}}*G_2\underset{G_{2,3}}*\cdots\underset{G_{j-2,j-1}}*G_{j-1}\underset{G_{j-1,j}}* G_j\underset{G_{j,j+1}} *\cdots\underset{G_{m-1,m}}*G_m
\end{eqnarray}

Then the chain of these building balls $\{B_j\}, 1\leq j\leq m$,
contains the bounded polyhedron $P_1\subset\Omega_1$, and the unbounded polyhedron $P_2\subset\Omega_2$ is inside of the chain of the balls $\{\widehat{B_j}, 1\leq j\leq m\}$, which are the complements in $S^3$ to the balls $B_j$.

\begin{figure}
\centering
\epsfxsize=10cm
\epsfbox{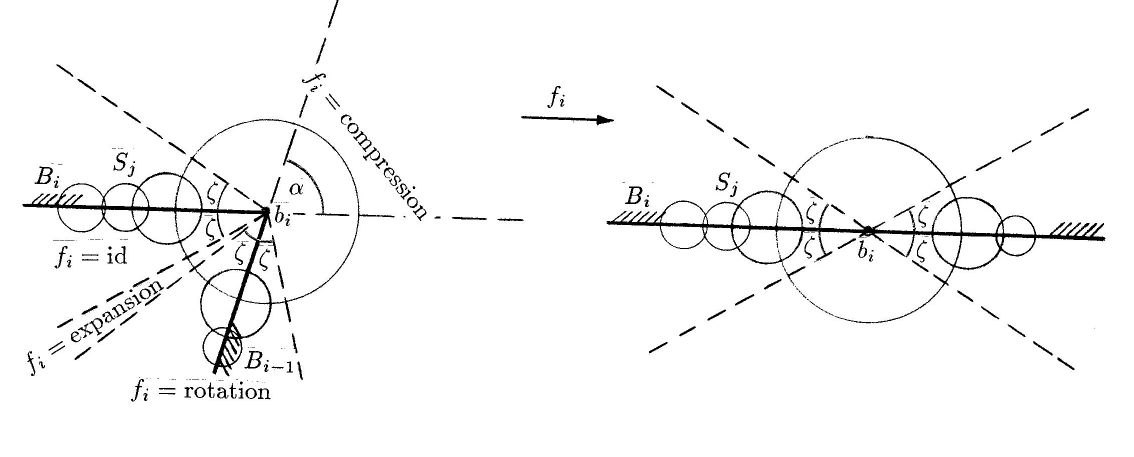}
\caption{Elementary bending homeomorphism $f_i$.}
\label{fig9}
\end{figure}

Now for each pair of balls $B_{i-1}$ and $B_i$
with the common boundary circle $b_i=\p B_{i-1}\cap\p B_i$, $1\leq i\leq m$, we construct a quasi-conformal
bending homeomorphism
$f_i$ that transfers $B_i\cup B_{i-1}$ onto the ball $B_i$ and which is
conformal in dihedral $\zeta_i$-neighborhoods of the spherical disks
$\p B_i\backslash\overline{B_{i-1}}$ and
$\p B_{i-1}\backslash\overline{B_i}$.
Namely, let $B_i$ and $B_{i-1}$ be half-spaces whose
boundary planes $\p B_i$ and $\p B_{i-1}$ contain the origin and intersect
along the third coordinate axis
$b_i=\{x\in\mathbb{R}^3\col x_1=x_2=0\}$ at an
angle $\alpha$, $0<\alpha<\pi$, and let $\zeta$ be a fixed number such
that $0< \zeta< \pi/2$ and
$0< \alpha< \pi - 2\zeta$.
Assuming the natural complex structure in
the plane $\mathbb{R}^2=\{x\in\mathbb{R}^3\col x_3=0\}$, we define
the quasi-conformal elementary bending homeomorphism $f_i$ by its
restriction to
the plane $\mathbb{C}=\mathbb{R}^2$ (see Figure \ref{fig9}), where

\begin{eqnarray}\label{bend}
f_i(z)=
\begin{cases}
z & \text{if}\ |\arg z|\geq\pi-\zeta\\
z\exp\left(i(\pi-\alpha)\right) & \text{if}\ |\pi-\alpha-\arg z|\leq\zeta\\
z\exp\left(i(\pi-\alpha)(1-\frac{\arg z-\zeta}{\pi-2\zeta})\right) & \text{if}\ \alpha-\pi+\zeta<\arg z<\pi-\zeta\\
z\exp\left(i(\pi-\alpha)(1+\frac{\zeta+\arg z}{\pi-2\zeta})\right) & \text{if}\ \zeta-\pi<\arg z<\alpha-\pi-\zeta
\end{cases}
\end{eqnarray}

Note that in each $i$-th step, $2\leq i\leq m$, we reduce the number of
balls $B_j$ in our chain by one. The constructed quasiconformal homeomorphisms $f_i$
have the properties:
\begin{enumerate}
\item  $f_i=\id$ in a neighborhood of reflection spheres from our collection $\Sigma$
that are disjoint from the circle $b_i$ and intersect some balls $B_j$, $j\geq i$.
\item  The composition
$f_{m+1}f_i f_{i-1}\cdots f_2 f_1$ transfers all spheres
from $\Sigma$ to spheres orthogonal to the boundary spheres of some balls
$B_j$, $i\leq j\leq m$, where all intersection angles between these spheres do not
change.
\end{enumerate}

Finally, renormalizing our last ball $B_m$ as the unit ball $B(0,1)$, we define our desired quasiconformal homeomorphism
$\phi_1\col P_1\rightarrow P_0$ as the restriction of the composition
$f_{m+1}f_mf_{m-1}\cdots f_2f_1$ of our bending
homeomorphisms $f_j$ on the fundamental polyhedron $P_1\subset\Omega_1$. Similarly (working with the balls $\widehat{B_j}$)
 we define the second quasiconformal homeomorphism $\phi_2\col P_2\rightarrow \widehat{P_0}$. Both mappings preserve the combinatorial structure of the polyhedra and their dihedral angles.

\section{Locally homeomorphic quasiregular mappings}

Now we apply results of the previous Section 3 to define our quasiregular mapping $F$ from $S^3\setminus S_*$ onto $S^3$.

\begin{theorem}\label{map}
 Let the uniform hyperbolic lattice $\Gamma\subset\is H^3$ and its discrete representation $\rho\col\Gamma\rightarrow G\subset\is H^4$ with the kernel as a free subgroup $F_3\subset\Gamma$ be as in Theorem \ref{constr}. Then there is a locally homeomorphic quasiregular mapping $F\col S^3\setminus S_*\rightarrow S^3$ whose all singularities lie in an exceptional subset $S_*$ of the unit sphere $S^2\subset S^3=\mathbb{R}^3 \cup\{\infty\}$ and form a dense in $S^2$ $\Ga$-orbit of a Cantor subset with Hausdorff dimension $\ln 5/\ln 6 \approx 0.89822444$. These (essential) singularities create a barrier for $F$ in the sense that at any such point $x\in S_*$ the map $F$ does not have radial limits on either side of $S^2\subset S^3$.
\end{theorem}

\begin{pf}
 First we define our locally homeomorphic quasiregular mapping $F$ in the complement $S^3\setminus S^2$ of the unit sphere,
 $F\col S^3\setminus S^2\rightarrow \Omega(G)\subset S^3$, $F(\infty)=\infty$.

We recall that the discrete group $G=\rho(\Gamma)\subset\is H^4\cong\Mob (3)$ constructed in Section 2 (Theorem \ref{constr}) has its discontinuity set $\Omega(G)\subset S^3$, and its conformal and discontinuous action in $\Omega(G)$ has $P=P_1\cup P_2$ as the fundamental polyhedron. Its symmetric connected components $P_1$ and  $P_2$ constructed in Section 2 have the combinatorial type of the convex 3-dimensional hyperbolic polyhedron $P_0$ fundamental for our hyperbolic lattice $\Gamma\subset\is H^3$ (for its conformal action in the unit ball $B^3(0,1)$). Let $\widehat{P_0}$ be the symmetric image of $P_0\subset B^3(0,1)$ with respect to the reflection in the unit sphere $S^2=\p B^3$. The polyhedron $P_0\cup\widehat{P_0}$ is a fundamental polyhedron (with two connected components which are convex in the induced hyperbolic metrics) for conformal and discontinuous action of our hyperbolic lattice $\Gamma$ in $S^3\setminus S^2$.

 In the previous Section 3 we have constructed quasiconformal homeomorphisms $\phi_1^{-1}\col P_0\rightarrow P_1$ and
$\phi_2^{-1}\col \widehat{P_0}\rightarrow P_2$.
 These two homeomorphisms map polyhedral sides of the polyhedra fundamental for the $\Gamma$-action to the corresponding sides of the polyhedra fundamental for the $G$-action and preserve combinatorial structures of polyhedra as well as their dihedral angles. Equivariantly extending these homeomorphisms, we define a quasiregular mapping
$F\col S^3\setminus S^2\rightarrow \Omega(G)$:

\begin{eqnarray}\label{QR}
F(x)=
\begin{cases}
\rho(\gamma)\circ \phi_1^{-1}\circ \gamma^{-1}(x)  & \text{if}\ |x|<1,\, x\in\gamma(P_0), \,\gamma\in\Gamma \\
\rho(\gamma)\circ \phi_2^{-1}\circ \gamma^{-1}(x)  & \text{if}\ |x|>1,\, x\in\gamma(\widehat{P_0}), \,\gamma\in\Gamma
\end{cases}
\end{eqnarray}

Since the initial quasiconformal homemorphisms $\phi_1^{-1}$ and $\phi_2^{-1}$ preserve combinatorial structures of polyhedra and their dihedral angles, the tesselations of $\Omega(\Gamma)=S^3\setminus S^2$ and $\Omega(G)\subset S^3$ by corresponding $\Gamma$- and $G$-images of fundamental polyhedra of the reflection groups $\Gamma$ and $G$ around all sides of polyhedra including their edges and vertices are perfectly similar. This implies that our quasiregular mapping $F$ defined by (\ref{QR}) is locally homeomorphic and $F(\infty)=\infty$.

It follows from Lemma \ref{h-body} that the limit set $\Lambda(G)\subset S^3$ of the group $G\subset \Mob(3)$ defines a Heegard splitting of infinite genus of the 3-sphere $S^3$ into two connected components $\Om_1$ and $\Om_2$ of the discontinuity set $\Om(G)$. The action of $G$  on the limit set $\La(G)$ is an ergodic word hyperbolic action (quasi-self-similar in the sense of Sullivan, see \cite{A2}, Cor. 2.66). For this ergodic action the set of fixed points of loxodromic elements $g\in G$ (conjugate to similarities in $\mathbb{R}^3$) is dense in $\La(G)$. Preimages $\ga\in\Ga$ of such loxodromic elements $g\in G$ for our homomorphism $\rho\col\Ga\rightarrow G$ are loxodromic elements in $\Ga$ with two fixed points $p,q\in\La(\Ga)=S^2, p\neq q$. This and Tukia's arguments of the group completion (see \cite{T1} and \cite{A2}, Section 4.6) show that our mapping $F$ can be continuously extended to the set of fixed points of such elements $\ga\in\Ga$, $F(Fix(\ga))=Fix(\rho(\ga))$. The sense of this continuous extension is that if $\ga\in\Ga$ is a loxodromic preimage of a loxodromic element $g\in G$, $\rho(\ga)=g$, and if $x\in S^3\backslash S^2$ tends to its fixed points $p$ or $q$ along the hyperbolic axis of $\ga$ (in $B(0,1)$ or in its complement $\widehat{B(0,1)}$) (i.e. radially) then $\lim_{|x|\to 1} F(x)$ exists and equals to the corresponding fixed point of the loxodromic element $g=\rho(\ga)\in G$. In that sense one can say that the limit set $\La(G)$ (the common boundary of the connected components $\Om_1, \Om_2 \subset\Om(G)$) is the $F$-image of points in the unit sphere $S^2\subset S^3$. So the mapping $F$ is onto the whole sphere $S^3$.

Nevertheless not all loxodromic elements $\ga\in\Ga$ in the hyperbolic lattice $\Ga\subset\is H^3$ have their images $\rho(\ga)\in G$ as loxodromic elements. Proposition \ref{homo} shows that $\ker \rho\cong F_3$ is a free subgroup on three generators in the lattice $\Ga$, and all elements $\ga\in F_3$ are loxodromic. Now we look at radial limits $\lim_{x\to p}F(x)$ when $x$ radially tends to a fixed point $p\in S^2$ of this loxodromic element $\ga\in F_3\subset\Ga$.

\begin{figure}
\centering
\includegraphics[width=5cm]{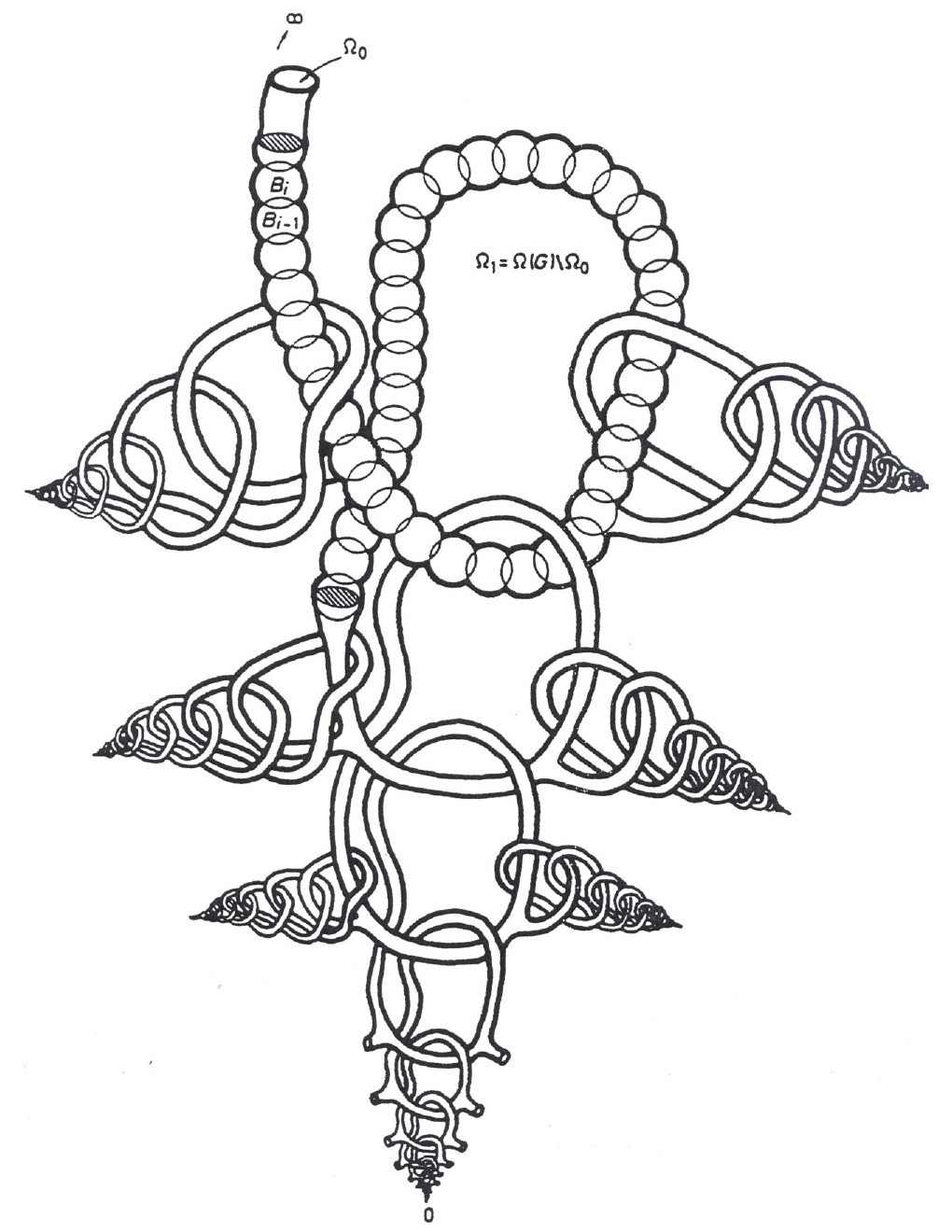}
\caption{Locally inextensible wild embedding of a closed ball into $\mathbb{R}^3$.}
\label{fig5}
\end{figure}
For a group $\Ga$ with a finite set $\Sigma=\{\ga_1, \ldots \ga_k\}$ of generators we consider its Cayley graph $K(\Ga,\Sigma)$, i.e. a 1-complex whose set of vertices is $\Ga$ and such that $a, b \in \Ga$ are joined by an edge if and only if
$a=b\ga^{\pm 1}$ for some $\ga\in\Sigma$. Since our $\Ga$ is a co-compact lattice acting in the hyperbolic space, we may define an embedding $\varphi$ of its Cayley graph $K(\Ga,\Sigma)$ in the hyperbolic space $H^3$ (model in the unit ball $B(0,1)$ or in its complement $\widehat{B(0,1)}$). For a point $0\in H^3$ not fixed by any $\ga\in\Ga\backslash\{1\}$, vertices $\ga\in K(\Ga,\Sa)$ are mapped to $\ga(0)$, and edges joining vertices
$a,b\in K(\Ga,\Sa)$ are mapped to the hyperbolic geodesic segments $[a(0),b(0)]$.
In other
words, $\varphi(K(\Ga,\Sa))$ is the graph that is dual to the tessellation of $H^3$ by
polyhedra $\ga(P_0)$ (or $\widehat{P_0}$), $\ga\in \Ga$.  Obviously, the map $\varphi$ is a $\Ga$-equivariant proper
embedding:  for any compact $C\subset H^3$, its pre-image
$\varphi^{-1}(\varphi(K(\Ga,\Sa))\cap C)$ is compact. Moreover this embedding is a pseudo-isometry (see \cite{C} and \cite{A2}, Theorem 4.35):

\begin{theorem}\label{embed}
 Let $\Ga\subset\is H^n$ be a convex co-compact group.
Then the map $\varphi\col K(\Ga,\Sa)\hookrightarrow H^n$ is a pseudo-isometry of the word metric
$(\ast,\ast)$ on $K(\Ga,\Sa)$ and the hyperbolic metric $d$, that is,
there are positive constants $K$ and $K'$ such that
\begin{eqnarray}\label{PI}
(a,b)/K\leq d(\varphi(a),\varphi(b))\leq K\cdot(a,b)
\end{eqnarray}
for all $a,b\in K(\Ga,\Sa)$ satisfying one of the following two conditions:
either $(a,b)\geq K'$ or $d(\varphi(a),\varphi(b))\geq K'$.
\end{theorem}

Theorem \ref{embed} implies (see \cite{A2}, Theorem 4.38) that the limit set of any convex-cocompact group
$\Ga\subset\Mob (n)$ can be identified with its group completion $\overline{\Ga}$, $\overline{\Ga}=\overline{K(\Ga,\Sa)}\setminus K(\Ga,\Sa)$. Namely there exists a continuous and $\Ga$-equivariant bijection $\varphi_{\Ga}\col \overline{\Ga} \rightarrow \La(\Ga)$.

Now for the kernel subgroup $F_3=\ker \rho \subset \Ga\subset\is H^3$ and for the pseudo-isometric embedding $\varphi$ from Theorem \ref{embed}, we consider its Cayley subgraph in $\varphi(K(\Ga,\Sa))\subset H^3$ which is a tree - see Figure \ref{fig10}. Since the limit set of $\ker \rho=F_3\subset\Ga$ corresponds to the 'bondary at infinity' $\p_{\infty} F_3$ of $F_3\subset\Ga$ (the group completion $\overline{F_3}$), it is a closed Cantor subset of the unit sphere $S^2$ with Hausdorff dimension $\ln 5/\ln 6 \sim 0.89822444$ (zero measure).

\begin{remark}
 The most common construction of the used Cantor set $\p_{\infty} F_3$ (similar to the original Cantor's ternary construction of a perfect set that is nowhere dense) follows from the Figure \ref{fig10}. It can be obtained by removing one middle sixth of a line segment and repeating this for remaining five subsegments. This gives its Hausdorff dimension $\ln 5/\ln 6$.
\end{remark}

 The $\Ga$-orbit $\Ga(\La(F_3))$ of our Cantor set is a dense subset $S_*$ of $S^2=\La(\Ga)$ because of density in the limit set $\La(\Ga)$ of the $\Ga$-orbit of any limit point. In particular we have such dense $\Ga$-orbit $\Ga(\{p,q\})$ of fixed points $p$ and $q$ of a loxodromic element $\ga\in F_3\subset \Ga$ (the images of $p$ and $q$ are fixed points of $\Ga$-conjugates of such loxodromic elements $\ga\in F_3\subset \Ga$).

On the other hand let $x\in l_{\ga}$ where $l_{\ga}$ is the hyperbolic axis of an element $\ga\in F_3\subset\Ga$ (or its $\Ga$-conjugate). Here the hyperbolic axes are either in $B(0,1)$ or in its complement $\widehat{B(0,1)}$). Denoting $d_{\ga}$ the translation distance of $\ga$, we have that any segment $[x, \ga(x)]\subset l_{\ga}$ is mapped by our quasiregular mapping $F$ to a non-trivial closed loop $F([x, \ga(x)])\subset\Om(G)=\Om_1\cup\Om_2$, inside of a handle of the mutually linked handlebodies $\Om_1$ or $\Om_2$ (similar to the loops $\beta_1\subset \Om_1$ and $\beta_2\subset \Om_2$ constructed in the proof of Lemma \ref{h-body}). Therefore when $x\in l_{\ga}$ radially tends to a fixed point $p$ (in $\fix (\ga)\in S^2$) of such element $\ga$, its image $F(x)$ goes along that closed loop $F([x, \ga(x)])\subset\Om(G)$ because $F(\ga(x))=\rho(\ga)(F(x))=F(x)$. Immediately it implies that the radial limit $\lim_{x\to p} F(x)$ does not exist. This shows that fixed points of any element $\ga\in F_3\subset \Ga$ (or its conjugate) are essential (topological) singularities of our quasiregular mapping $F$. So our quasiregular mapping $F$ has no continuous extension to the subset $S_*\subset S^2$ (from both sides of the unit sphere $S^2$ in $S^3$) which is a dense subset of the unit sphere $S^2\subset S^3$.
\end{pf}

\begin{remark}
In terms of the holomorphic function theory of several complex variables, both components of the complement $S^3\setminus S^2$ play the role of the so called domain of holomorphy for the constructed in Theorem \ref{map} locally homeomorphic quasiregular mapping $F$. Obviously instead of the sphere $S^2\subset S^3$ one may consider any quasi-sphere $S^2_q\subset S^3$ which is the image of $S^2$ under a quasiconformal homeomorphism of $S^3$. In other words, the complement $S^3\setminus S^2_q$ consisting of two quasi-balls has the same property of domain of holomorphy for a locally homeomorphic quasiregular mapping $F$ constructed in Theorem \ref{map}, and the quasi-sphere $S^2_q\subset S^3$ is a barrier for it with a dense subset of essential singularities.
\end{remark}

\begin{remark}
The constructed in Theorem \ref{map} barrier $S^2$ for our locally homeomorphic quasiregular mapping $F$ in the 3-space has completely different nature from the topological barrier $S^2$ for the quasisymmetric embedding $f: \overline{B(0,1)}\hookrightarrow \mathbb{R}^3$ of the closed unit ball $\overline{B(0,1)}$ into $\mathbb{R}^3$ constructed in \cite{A3}. That topological barrier for the embedding $f$ was due to wild knottings of the boundary topological sphere $f(S^2)\subset \mathbb{R}^3$ on its dense subset. Figure \ref{fig5} demonstrates the topological nature of this wild knotting when the fundamental group of the complement
$\mathbb{R}^3\setminus f(\overline{B(0,1)})$ is infinitely generated near wild knotting points.
\end{remark}

  Dept of Math., University of Oklahoma, Norman, OK 73019, USA

   e-mail: apanasov\char`\@ ou.edu


\vskip20pt


\begin{thebibliography}{99}
\bibitem[1]{An1} E.M.~Andreev, \textit{On convex polyhedra in Lobachevsky space}, Mat. Sbornik
\textbf{81}, 1970, 445--478 (In Russian); Engl. Transl.: Math. USSR Sbornik
\textbf{10}, 1970, 413--440.

\bibitem[2]{An2}  E.M.~Andreev,
\textit{The intersections of the planes of faces of polyhedra with
sharp angles}, Mat. Zametki \textbf{8}, 1970, 521--527 (in Russian);
Engl. Transl.: Math. Notes \textbf{8}, 1970, 761--764.


\bibitem[3]{A1}  Boris Apanasov,
\textit{Nontriviality  of Teichm\"uller space for
Kleinian  group in  space.} - In: Riemann Surfaces and Related
Topics, Proc. 1978 Stony  Brook Conference (I.Kra and
B.Maskit, eds), Ann. of  Math. Studies \textbf{97},
Princeton Univ. Press, 1981, 21--31.

\bibitem[4]{A2}  Boris Apanasov,
\textit{Conformal geometry of
discrete groups and manifolds.} - De Gruyter Exp. Math. \textbf{32}, W. de Gruyter, Berlin - New York, 2000.

\bibitem[5]{A3} Boris Apanasov,
\textit{Quasisymmetric embeddings of a closed ball
inextensible in neighborhoods of any boundary points}, Ann. Acad.
    Sci. Fenn., Ser. A I Math \textbf{14}, 1989, 243--255.

\bibitem[6]{A4} Boris Apanasov,
\textit{Nonstandard uniformized conformal
structures on hyperbolic manifolds.} - Invent. Math. \textbf{105}, 1991, 137--152.

\bibitem[7]{A5} Boris Apanasov,
\textit{Hyperbolic 4-cobordisms and group homomorphisms with infinite kernel.}
- Atti Semin. Mat. Fis. Univ. Modena Reggio Emilia \textbf{57}, 2010, 31--44.

\bibitem[8]{AT} Boris N. Apanasov and Andrei V. Tetenov,
\textit{Nontrivial cobordisms with
    geometrically finite hyperbolic structures}, J. of Diff. Geom. \textbf{28},
 1988, 407--422.

\bibitem[9]{B} Pavel P.~Belinskii
\textit{On the continuity of quasiconformal mappings in space and Liouville's theorem}, Dokl. Akad. Nauk SSSR
 \textbf{147}, 1962, 1003--1004. (in Russian).

\bibitem[10 ]{BL1} Pavel P.~Belinskii and Mikhael A.~Lavrentiev,
\textit {Certain problems of the geometric function theory},
Trudy Math. Inst. Steklov \textbf{128:2}, Nauka, Moscow, 1972, 34--40 (in Russian);
Engl. Translation: Proc. Steklov Inst. Math. \textbf{128}, Amer. Math. Soc., Providence, 1973.

\bibitem[11]{BL2} Pavel P.~Belinskii and Mikhael A.~Lavrentiev,
\textit {On locally quasiconformal mappings in $(n\geq 3)$-space}, In:
Contributions to Analysis: a collection of papers dedicated to Lipman Bers
 (L.~Ahlfors et al, eds.), Academic Press, 1974, 27--30.

\bibitem[12] {C} James W.~Cannon, \textit{The combinatorial structure of cocompact discrete
hyperbolic groups}, Geom. Dedicata \textbf{16}, 1984, 123--148.

\bibitem[13]{DP} David Drasin and Pekka Pankka, \textit{Sharpness of Rickman’s Picard
theorem in all dimensions}, Acta Math. \textbf{214}, 2015, 209--306.

\bibitem[14]{F} P. Fatou,
\textit{S\'{e}ries trigonom\'{e}triques et s\'{e}ries de Taylor}, Acta Math. \textbf{30}, 1906, 335--400 (in French).

\bibitem[15]{Go} Vladimir M.~Goldstein,
\textit{The behavior of mappings with bounded distortion when the distortion coefficient is close to one}, Sibirsk. Mat. Z. 12, 1971, 1250--1258 (in Russian).

\bibitem[16]{G1} Mikhael Gromov,
\textit {Hyperbolic manifolds, groups and actions.} -  In: Riemann Surfaces and Related
Topics: Proc. 1978 Stony  Brook Conference (I.Kra and B.Maskit, eds.), Ann. of  Math. Studies \textbf{97},
Princeton Univ. Press, 1981, 183--213.

\bibitem[17]{G2} Mikhael Gromov,
\textit {Structures m\'{e}triques pour les vari\'{e}t\'{e}s riemanniennes.} - 
Cedic, Paris, 1981.

\bibitem[18]{G3} Mikhael Gromov,  
\textit{Metric structures for Riemannian and non-Riemannian spaces}, Birkh\"{a}user, Boston, 1999.

\bibitem[19]{Gr} H.Gr\"{o}zsch,
\textit{\"{U}ber die Verzerrung bei schlichten nichtkonformen Abbildungen und \"{u}ber eine damit
zusammenh\"{a}ngende Erweiterung des Picardischen Satzes},
Ber. Verh. S\"{a}chs.Acad. Wiss. Leipzig \textbf{80}, 1928, 503--507.

\bibitem[20]{L1} Mikhael A.~Lavrentiev,
\textit{On a class of continuous mappings}, Mat. Sbornik, \textbf{42}, 1935, 407--424. (in Russian).

\bibitem[21]{L2} Mikhael A.~Lavrentiev,
\textit{On a differential test for homeomorphism of mappings of three-dimensional
domains}, Doklady Akad. Nauk SSSR \textbf{20}, 1938, 241--242 (in Russian).

\bibitem[22]{MR} Olli Martio and Seppo Rickman,
\textit{Boundary behavior of quasiregular mappings}, Ann. Acad. Sci. Fenn. Ser. A I Math
\textbf{507}, 1972, 1--17.

\bibitem[23]{MRV} Olli Martio, Seppo Rickman and Jussi V\"{a}is\"{a}l\"{a},
\textit {Distortion and singularities of quasiregular mappings},
Ann. Acad. Sci. Fenn. Ser. A I Math \textbf{465}, 1970, 1--13.

\bibitem[24]{MS1} Olli Martio and Uri Srebro,
\textit{Automorphic quasimeromorphic mappings in $\mathbb{R}^n$},   .
Acta Math. \textbf{135}, 1975, 221--247.

\bibitem[25]{MS2} Olli Martio and Uri Srebro,
\textit{Locally injective automorphic mappings in $\mathbb{R}^n$},
Math. Scand. \textbf{85}, 1999, 49--70.

\bibitem[26]{Re} Yuri G. Reshetnyak,
\textit {Space mappings with bounded distortion.} - Nauka, Novosibirsk, 1982 (in Russian); Engl. transl.:
Transl. Math. Monogr. \textbf{73}, Amer. Math. Soc., Providence, R.I., 1989.

\bibitem[27]{Ra} Kai Rajala,
\textit{Radial limits of quasiregular local homeomorphisms}, Amer. J. Math. \textbf{130}, 2008, 269--289.

\bibitem[28]{RS} Luis Ribes and Benjamin Steinberg,
\textit{A wreath product approach to classical subgroup theorems},
Enseign. Math. (2) \textbf{56}, 2010),  49--72.

\bibitem[29]{Ri} Seppo Rickman,
\textit {Quasiregular Mappings.} - Ergeb. Math. Grenzgeb. \textbf{26},
Springer, Berlin–Heidelberg, 1993.

\bibitem[30]{S}  Dennis Sullivan,
\textit {Quasiconformal homeomorphisms and dynamics, II: Structural
stability implies hyperbolicity for Kleinian groups}, Acta Math. \textbf{155}, 1985, 243--260.

\bibitem[31]{T1} Pekka Tukia,
\textit{On isomorphisms of geometrically Kleinian groups}, Publ. Math. IHES \textbf{61}, 1985,
171--214.

\bibitem[32]{T2} Pekka Tukia,
\textit{Automorphic quasimeromorphic mappings for torsionless hyperbolic groups},
Ann. Acad. Sci. Fenn. Ser. A I Math. \textbf{10}, 1985, 545--560.


\bibitem[33]{V} Jussi V\"{a}is\"{a}l\"{a},
\textit{A survey of quasiregular maps in $\mathbb{R}^n$}, In:
Proc. ICM (Helsinki, 1978), Acad. Sci. Fennica, Helsinki, 1980, 685--691.

\bibitem[34]{Vu} Matti Vuorinen,
\textit{Conformal geometry and quasiregular mappings.} - Lecture Notes in Math \textbf{1319}, Springer,
 Berlin–Heidelberg, 1988.

\bibitem[35]{Z1} Vladimir A.~Zorich,
\textit{A theorem of M. A. Lavrentiev on quasiconformal space maps}, Mat. Sbornik \textbf{74},
1967, 417--433 (in Russian); Engl. transl.: Math. USSR Sbornik, \textbf{3} (1967), 389--403.

\bibitem[36]{Z2} Vladimir A.~Zorich,
\textit{On Gromov’s geometric version of the global homeomorphism theorem for
quasi-conformal mappings.} - Abstracts, XIV Rolf Nevanlinna Colloquium, Helsinki, June
10--14, 1990, 36.

\bibitem[37]{Z3} Vladimir A.~Zorich,
\textit{Quasi-conformal maps and the asymptotic geometry of manifolds}, Uspekhi Mat. Nauk \textbf{57:3}, 2002, 3–-28 (in Russian); Engl. transl.: Russian Math. Surveys \textbf{71}, 2016, 161–-163.

\bibitem[38]{Z4} Vladimir A.~Zorich,
\textit{A note on the injectivity radius of quasiconformal immersions}, Uspekhi Mat. Nauk \textbf{71:1}, 2016, 173--174 (in Russian); Engl. transl.: Russian Math. Surveys \textbf{71}, 2016, 161–-163.

\bibitem[39]{Z5} Vladimir A.~Zorich,
\textit{Some observations concerning multidimensional quasiconformal mappings}, Mat. Sbornik, \textbf{208:3}, 2017, 72–-95
(in Russian); Engl. transl.: Sb. Math. \textbf{208}, 2017, 377--398.


\end{thebibliography}
\end{document}